\font\fivemsb=msbm5
\font\sevenmsb=msbm7
\font\tenmsb=msbm10
\newfam\msbfam
\textfont\msbfam=\tenmsb
\scriptfont\msbfam\sevenmsb
\scriptscriptfont\msbfam\fivemsb
\def\sqbox{{\vcenter{\hrule height .5pt \hbox{\vrule
  width .5pt height 5pt \kern 5pt \vrule width .5pt}
  \hrule height .5pt}}}
\def\qed{\sqbox}

\def\a{\alpha}
\def\b{\beta}
\def\d{\delta}
\def\cv{{\cal V}}

\documentstyle[verbatim]{amsart}

\theoremstyle{plain}
\newtheorem{theorem}{Theorem}

\newtheorem{lemma}{Lemma}
\newtheorem{bctheorem}{Bender-Canfield Theorem \cite{bc}}

\theoremstyle{definition}

\renewcommand{\rm}{\normalshape}

\begin{document}
\title[Bell numbers]%
   {Bell numbers, log-concavity, and log-convexity}

\author{Nobuhiro Asai}
\author{Izumi Kubo}
\author{Hui-Hsiung Kuo}
\address{Nobuhiro Asai: Graduate School of Mathematics\\
  Nagoya University \\ Nagoya 464-8602 \\ JAPAN}
\address{Izumi Kubo: Department of Mathematics \\
Faculty of Science \\ Hiroshima University \\
Higashi-Hiroshima 739-8526 \\ JAPAN}
\address{Hui-Hsiung Kuo: Department of Mathematics\\
  Louisiana State University \\ Baton Rouge \\
LA 70803 \\ USA}

\thanks{Research supported by the Daiko Foundation 1998
(N.A.), U.S.~Army Research Office grant \#DAAH04-94-G-0249,
Academic Frontier in Science of Meijo University, and
the National Science Council of Taiwan (H.-H.K.)}

\maketitle

\begin{abstract}
Let $\{b_{k}(n)\}_{n=0}^{\infty}$ be the Bell numbers
of order $k$. It is proved that the sequence
$\{b_{k}(n)/n!\}_{n=0}^{\infty}$ is log-concave
and the sequence $\{b_{k}(n)\}_{n=0}^{\infty}$ is
log-convex, or equivalently, the following inequalities
hold for all $n\geq 0$,
$$1\leq {b_{k}(n+2) b_{k}(n) \over b_{k}(n+1)^{2}} \leq
{n+2 \over n+1}.$$
Let $\{\a(n)\}_{n=0}^{\infty}$ be a sequence of positive
numbers with $\a(0)=1$. We show that if
$\{\a(n)\}_{n=0}^{\infty}$ is log-convex, then
$$\a (n) \a (m) \leq \a(n+m), \quad \forall n, m\geq 0.$$
On the other hand, if $\{\a(n)/n!\}_{n=0}^{\infty}$ is
log-concave, then
$$\a (n+m) \leq {n+m \choose n} \a (n) \a (m), \quad
\forall n, m\geq 0.$$
In particular, we have the following inequalities for the
Bell numbers
$$b_{k}(n) b_{k}(m) \leq b_{k}(n+m) \leq {n+m \choose n}
  b_{k}(n) b_{k}(m),  \quad \forall n, m\geq 0.$$
Then we apply these results to white noise distribution
theory.
\end{abstract}

\smallskip
\section{The main theorems} \label{sec:1}

For an integer $k\geq 2$, let $\exp_{k}(x)$ denote the
$k$-times iterated exponential function
$$\exp_{k} (x) = \underbrace{\exp\bigl(\exp\cdots\bigl(
\exp(x) \bigr)\bigr)}_{k-times}.$$
Let $\{B_{k}(n)\}_{n=0}^{\infty}$ be the sequence of numbers
given in the power series of $\exp_{k}(x)$
\begin{equation} \label{eq:n-9}
\exp_{k} (x) = \sum_{n=0}^{\infty} {B_k(n)
\over n!}\,x^{n}.
\end{equation}
The {\em Bell numbers $\{b_{k}(n)\}_{n=0}^{\infty}\,$
of order $k\,$} are defined by
$$b_k(n) = {B_{k}(n) \over \exp_{k}(0)}, \quad  n\geq 0.$$
The numbers $b_{2}(n), n\geq 0,$ with $k=2$ are usually
known as the {\em Bell numbers}. The first few terms of
these numbers are $1, 1, 2, 5, 15, 52, 203$. Note that
$\exp_{2}(0)=e$ and so we have
\begin{equation} \label{eq:1-1}
e^{e^x-1}=\sum_{n=0}^{\infty} {b_2(n) \over n!}\,x^{n}.
\end{equation}

A sequence $\{\d(n)\}_{n=0}^{\infty}$ of nonnegative real
numbers is called {\em log-concave} if
$$\d(n) \d(n+2) \leq \d(n+1)^{2}, \quad \forall n\geq 0.$$
It is called {\em log-convex} if
$$\d(n) \d(n+2) \geq \d(n+1)^{2}, \quad \forall n\geq 0.$$
\par
The main purpose of this paper is to prove the following
theorems.

\begin{theorem} \label{thm:1}
Let $\{b_{k}(n)\}_{n=0}^{\infty}$ be the Bell numbers of
order $k$. Then the sequence $\{b_{k}(n)/n!\}_{n=0}^{\infty}$
is log-concave and the sequence $\{b_{k}(n)\}_{n=0}^{\infty}$
is log-convex.
\end{theorem}

Note that the conclusion of this theorem is equivalent the
inequalities
\begin{equation}
1\leq {b_{k}(n) b_{k}(n+2) \over b_{k}(n+1)^{2}} \leq
{n+2 \over n+1}, \quad \forall n\geq 0.  \notag
\end{equation}

A different proof of the log-convexity of
$\{b_{2}(n)\}_{n=0}^{\infty}$ has been given earlier by
Engel \cite{engel}. In \cite{can} Canfield showed that
the log-concavity of $\{b_{2}(n)/n!\}_{n=0}^{\infty}$
holds asymptotically. In a recent paper \cite{cks}, Cochran
et al.~used the log-concavity of certain sequences to
study characterization theorems. However, they did not show
whether the sequence $\{b_{k}(n)/n!\}_{n=0}^{\infty}$ is
log-concave. Thus our Theorem \ref{thm:1} fills up this gap
(for details, see Section \ref{sec:3}.)

\begin{theorem} \label{thm:2}
Let $\{\a(n)\}_{n=0}^{\infty}$ be a sequence of positive
numbers with $\a(0)=1$.
\par\noindent
{\rm (a)} If $\{\a(n)\}_{n=0}^{\infty}$ is log-convex,
then
\begin{equation} \label{eq:n-7}
\a (n) \a (m) \leq \a(n+m), \quad \forall n, m\geq 0.
\end{equation}
{\rm (b)} If $\{\a(n)/n!\}_{n=0}^{\infty}$ is log-concave,
then
\begin{equation} \label{eq:n-8}
\a (n+m) \leq {n+m \choose n} \a (n) \a (m), \quad
\forall n, m\geq 0.
\end{equation}
\end{theorem}

We will prove Theorems \ref{thm:1} and \ref{thm:2} in
Section \ref{sec:2}. The next theorem is an immediate
consequence of these two theorems.

\begin{theorem} \label{thm:3}
The Bell numbers $\{b_{k}(n)\}_{n=0}^{\infty}$ of order
$k$ satisfy the inequalities
\begin{equation} \label{eq:n-2}
b_{k}(n) b_{k}(m) \leq b_{k}(n+m) \leq {n+m \choose n}
  b_{k}(n) b_{k}(m),  \quad \forall n, m\geq 0.
\end{equation}
\end{theorem}

In a recent paper \cite{kks} it is shown that for any $k
\geq 2$ there exist constants $c_{2}$ and $c_{3}$,
depending on $k$, such that for all $n, m\geq 0$,
$$b_{k}(n+m) \leq c_{2}^{n+m} b_{k}(n)  b_{k}(m),
\qquad  b_{k}(n)  b_{k}(m) \leq c_{3}^{n+m} b_{k}(n+m).$$
\par
Observe that from Eq.~(\ref{eq:n-2}) we get $b_{k}(n)
b_{k}(m) \leq b_{k}(n+m) \leq 2^{n+m} b_{k}(n) b_{k}(m)$.
Thus in fact we can take $c_{2}=2$ and $c_{3}=1$ for
the Bell numbers of any order $k$.

\smallskip
\section{Proofs of Theorems \ref{thm:1} and
\ref{thm:2}} \label{sec:2}

For the proof of Theorem \ref{thm:1} we prepare two lemmas
and state the Bender-Canfield theorem \cite{bc}.

\begin{lemma} \label{lem:1}
If $\{\b(n)/n!\}_{n=0}^{\infty}$ is a log-concave sequence
and $r$ is a nonnegative real number such that $\b(2) \leq
r \b(1)^{2}$, then the sequence $1, \, r\b(n)/(n-1)!, \>
n\geq 1,$ is log-concave.
\end{lemma}

\begin{pf}
By assumption we have
\begin{equation}
{\b(n) \over n!} {\b(n+2) \over (n+2)!} \leq
\bigg({\b(n+1) \over (n+1)!}\bigg)^{2}.  \notag
\end{equation}
When $n\geq 1$ this inequality is equivalent to
\begin{equation}
{(n+1)^{2} \over n(n+2)}\bigg({\b(n) \over (n-1)!}
{\b(n+2) \over (n+1)!} \bigg) \leq
\bigg({\b(n+1) \over n!}\bigg)^{2}.  \notag
\end{equation}
Note that $(n+1)^{2} \geq n(n+2)$. Hence for $n\geq 1$,
\begin{equation}
  {\b(n) \over (n-1)!} {\b(n+2) \over (n+1)!} \leq
\bigg({\b(n+1) \over n!}\bigg)^{2}.  \notag
\end{equation}
Thus for any constant $r$ we have
\begin{equation}
{r \b(n) \over (n-1)!} {r \b(n+2) \over (n+1)!} \leq
\bigg({r \b(n+1) \over n!}\bigg)^{2}, \quad
\forall n \geq 1.  \notag
\end{equation}
Moreover, the assumption $\,\b(2) \leq r \b(1)^{2}$ implies
that $1 \!\cdot \!(r \b(2))\leq\big(r \b (1)\big)^{2}$. Thus
the sequence $1, \, r \b(n)/(n-1)!,\> n\geq 1,$ is
log-concave.
\end{pf}

\begin{bctheorem}
Let $1, Z_1, Z_2, \dots$ be a log-concave sequence of
nonnegative real numbers and define the sequence
$\{a(n)\}_{n=0}^{\infty}$ by
\begin{equation}
\sum_{n=0}^{\infty} {a(n) \over n!}\,x^{n} = \exp\bigg(
\sum_{j=1}^{\infty} {Z_{j} \over j}\,x^{j} \bigg). \notag
\end{equation}
Then the sequence $\{a(n)/n!\}_{n=0}^{\infty}$ is
log-concave and the sequence $\{a(n)\}_{n=0}^{\infty}$ is
log-convex.
\end{bctheorem}

\begin{lemma} \label{lem:2}
The sequence $\{b_{2}(n)/n!\}_{n=0}^{\infty}$ is log-concave
and the sequence $\{b_{2}(n)\}_{n=0}^{\infty}$ is log-convex.
\end{lemma}

\begin{pf}
Note that $e^{e^x-1} = \exp\big(\sum_{j=1}^{\infty}
{1 \over j!}\,x^{j}\big)$. Hence by Eq.~(\ref{eq:1-1})
we have
\begin{equation}
\exp\bigg(\sum_{j=1}^{\infty} {1 \over j!}\, x^{j}\bigg)
= \sum_{n=0}^{\infty} {b_2(n) \over n!}\,x^{n}.  \notag
\end{equation}
Let $Z_{j}={1\over (j-1)!}$ for $j\geq 1$. It is easy to
check that the sequence $1, Z_{1}, Z_{2}, \ldots$ is
log-concave. Thus this lemma follows from the above
Bender-Canfield theorem.
\end{pf}

\par\noindent
{\bf Proof of Theorem \ref{thm:1}}
\par
\smallskip
We prove the theorem by mathematical induction. By Lemma
\ref{lem:2} the theorem is true for $k=2$. Assume the
theorem is true for $k$. Note that
$\exp_{k+1}(x)/\exp_{k+1}(0) = \exp\big(\exp_{k}(x) -
\exp_{k}(0)\big)$. Hence
\begin{equation}
\exp\big(\exp_{k}(x) - \exp_{k}(0)\big) =
\sum_{n=0}^{\infty} {b_{k+1}(n) \over n!}\,x^{n}. \notag
\end{equation}
But $\exp_{k}(x) - \exp_{k}(0) = \sum_{j=1}^{\infty}
{B_{k}(j) \over j!}\, x^{j}$. Thus we get
\begin{equation}
\exp\bigg(\sum_{j=1}^{\infty} {B_{k}(j) \over j!}\, x^{j}
\bigg) = \sum_{n=0}^{\infty} {b_{k+1}(n) \over n!} \,
x^{n}.  \notag
\end{equation}
Let $Z_{j} = {B_{k}(j) \over (j-1)!},\, j\geq 1$.
Then the above equation becomes
\begin{equation} \label{eq:2-1}
\exp\bigg(\sum_{j=1}^{\infty} {Z_{j} \over j}\,x^{j}\bigg)
= \sum_{n=0}^{\infty} {b_{k+1}(n) \over n!}\,x^{n}.
\end{equation}
By the induction assumption, the sequence
$\{b_{k}(n)/n!\}_{n=0}^{\infty}$ is log-concave. This
implies in particular that $b_{k}(0) b_{k}(2)/2 \leq
b_{k}(1)^{2}$. But $b_{k}(0)=1$ and $\exp_{k}(0)>2$. Hence
$$b_{k}(2) \leq 2 b_{k}(1)^{2} < \exp_{k}(0)\,b_{k}(1)^{2}.$$
Thus we can apply Lemma \ref{lem:1} with $\b(n)=b_{k}(n)$ and
$r=\exp_{k}(0)$ to conclude that the sequence
$$1, \> \> \exp_{k}(0) {b_{k}(n) \over (n-1)!},\> n\geq 1,$$
is log-concave. Note that for $n\geq 1$,
$$\exp_{k}(0) {b_{k}(n) \over (n-1)!} = {B_{k}(n) \over
(n-1)!} = Z_{n}.$$
Hence the sequence $1, Z_{1}, Z_{2}, \ldots$ is log-concave.
Upon applying the Bender-Canfield theorem, we see from
Eq.~(\ref{eq:2-1}) that the sequence
$\{b_{k+1}(n)/n!\}_{n=0}^{\infty}$ is log-concave and
the sequence $\{b_{k+1}(n)\}_{n=0}^{\infty}$ is
log-convex.  \qed

\medskip\noindent
{\bf Proof of Theorem \ref{thm:2}}
\par
\smallskip
To prove (a), let $\{\a(n)\}_{n=0}^{\infty}$ be log-convex.
Then $\a(n) \a(n+2)\geq \a(n+1)^{2}$. Hence $\a(n+1)/\a(n)
\leq \a(n+2)/\a(n+1)$ and this implies that
for any $n\geq 0$ and $m \geq 1,$
$${\a(1) \over \a(0)} \leq {\a(2) \over \a(1)} \leq
\cdots \leq {\a(n+m) \over \a(n+m-1)}.$$
Therefore, for any $n\geq 0$ and $m \geq 1,$
$${\a(1) \over \a(0)} {\a(2) \over \a(1)} \cdots {\a(m)
\over \a(m-1)} \leq {\a(n+1) \over \a(n)} {\a(n+2)
\over \a(n+1)} \cdots {\a(n+m) \over \a(n+m-1)}.$$
After the cancellation we get $\a(n) \a(m) \leq \a(0)
\a(n+m)$. But $\a(0)=1$ and so Eq.~(\ref{eq:n-7}) is true
when $n\geq 0$ and $m \geq 1$. When $m=0$,
Eq.~(\ref{eq:n-7}) obviously holds for any $n\geq 0$.
Hence we have proved assertion (a).
\par
For the proof of (b), first note that
$\{\a(n)/n!\}_{n=0}^{\infty}$ is log-concave if and only
if for all $n\geq 0$,
$${\a(n+1) \over \a(n)} \geq
{n+1 \over n+2} {\a(n+2) \over \a(n+1)}.$$
By using this inequality repeatedly, we get the following
inequalities for any $n\geq 0$ and $m\geq 1$,
$${\a(1) \over \a(0)}
  \geq {1\over 2} {\a(2) \over \a(1)}
  \geq {1\over 3} {\a(3) \over \a(2)}
  \geq  \cdots
  \geq {1\over n+m} {\a(n+m) \over \a(n+m-1)}.$$
Hence for any $0 \leq j \leq m-1$,
$${\a(j+1) \over \a(j)} \geq {j+1 \over n+m}
{\a(n+m) \over \a(n+m-1)}.$$
Therefore,
\begin{align}
& {\a(1) \over \a(0)} {\a(2) \over \a(1)} \cdots
{\a(m) \over \a(m-1)}  \notag \\
& \geq
\bigg({1 \over n+1} {\a(n+1) \over \a(n)}\bigg)
\bigg({2 \over n+2} {\a(n+2) \over \a(n+1)}\bigg) \cdots
\bigg({m \over n+m} {\a(n+m) \over \a(n+m-1)}\bigg). \notag
\end{align}
After the cancellation we get
$${\a(m) \over \a(0)} \geq {n!m! \over (n+m)!}
{\a(n+m) \over \a(n)}.$$
But $\a(0)=1$. Hence we have proved that for any $n\geq 0,
m\geq 1$,
$$\a(n+m) \leq {n+m \choose n} \a(n) \a(m).$$
Note that when $m=0$, Eq.~(\ref{eq:n-8}) obviously holds
for any $n\geq 0$. Thus assertion (b) is proved.   \qed

\smallskip
\section{Application to white noise analysis} \label{sec:3}

\par\noindent
$\bullet$ {\em Characterization of test and generalized
functions}
\par
\smallskip
The Bell numbers $\{b_{k}(n)\}_{n=0}^{\infty}$ for $k\geq 2$
provide important examples in white noise distribution
theory \cite{kuo96}. In a recent paper \cite{cks} Cochran
et al.~have constructed a space $[\cv]_{\a}$ of test
functions and its dual space $[\cv]_{\a}^{*}$ of
generalized functions from a nuclear space $\cv$ and a
sequence $\{\a(n)\}_{n=0}^{\infty}$ of positive numbers
satisfying the following conditions:
\par
\smallskip
\begin{enumerate}
\item[(1)] $\a(0)=1$.
\smallskip
\item[(2)] $\inf_{n\geq 0} \a(n) >0$.
\smallskip
\item[(3)] $\lim_{n\to\infty} \big({\a(n) \over
n!}\big)^{1/n}\,=0$.
\end{enumerate}
\par
\smallskip
For the characterization of generalized functions in
$[\cv]_{\a}^{*}$ (Theorem 6.4 in \cite{cks}) they assume
the following condition
\begin{equation} \label{eq:1-2}
\limsup_{n\to\infty} \bigg({n! \over \a(n)} \,\inf_{x>0}
{G_{\a}(x) \over x^{n}}\bigg)^{1/n} < \infty.
\end{equation}
where $G_{\a}(x) = \sum_{n=0}^{\infty} {\a(n) \over n!}\,
x^{n}$ is the exponential generating function of the
sequence $\{\a(n)\}_{n=0}^{\infty}$. Furthermore,
by Corollary 4.4 in \cite{cks}, if the sequence
$\{\a(n)/n!\}_{n=0}^{\infty}$ is log-concave, then the
condition in Eq.~(\ref{eq:1-2}) is satisfied.
\par
For the case $\a(n) = b_{k}(n)$, Cochran et al.~showed
in Proposition 7.4 in \cite{cks} that the condition in
Eq.~(\ref{eq:1-2}) is satisfied. However, they did not
show whether the sequence $\{b_{k}(n)/n!\}_{n=0}^{\infty}$
is log-concave. Our Theorem \ref{thm:1} shows that this
is indeed the case.
\par
The other conclusion in Theorem \ref{thm:1}, i.e.,
$\{b_{k}(n)\}_{n=0}^{\infty}$ being log-convex,
can be used to characterize the test functions. First
we point out the following fact which can be easily
checked.
\par
\medskip\noindent
{\bf Fact.} {\em If $\{\b(n)\}_{n=0}^{\infty}$ is
log-convex, then $\big\{{1\over \b(n) n!}
\big\}_{n=0}^{\infty}$ is log-concave.}
\par
\medskip
Recall from Theorem \ref{thm:1} that the sequence
$\{b_{k}(n)\}_{n=0}^{\infty}$ is log-convex. Hence by
the above fact the sequence $\big\{{1\over b_{k}(n) n!}
\big\}_{n=0}^{\infty}$ is log-concave.
\par
In \cite{cks} Cochran et al.~did not study the
characterization of test functions in $[\cv]_{\a}$. In
a recent paper \cite{kubo} and our ongoing project
initiated in \cite{akk} several theorems on the
characterization of test functions and related results
have been obtained. For test functions, we need to
assume the following condition
\begin{equation} \label{eq:3-1}
\limsup_{n\to\infty} \bigg(\a(n) n! \,\inf_{x>0}
{G_{1/\a}(x) \over x^{n}}\bigg)^{1/n} < \infty,
\end{equation}
where $G_{1/\a}(x) = \sum_{n=0}^{\infty} {1 \over
\a(n) n!}\,x^{n}$ is the exponential generating function
of the sequence $\{{1 \over \a(n)}\}_{n=0}^{\infty}$.
The same argument as in the proof of Corollary 4.4 in
\cite{cks} can be used to show that if $\{{1\over
\a(n) n!}\}_{n=0}^{\infty}$ is log-concave, then the
condition in Eq.~(\ref{eq:3-1}) is satisfied.
\par
In particular, when $\a(n) = b_{k}(n)$, we know from
Theorem \ref{thm:1} that the sequence
$\{{1\over b_{k}(n) n!}\}_{n=0}^{\infty}$ is log-concave.
Thus the condition in Eq.~(\ref{eq:3-1}) is satisfied.
\par
\medskip\noindent
$\bullet$ {\em Inequality conditions on the sequence
$\{\a(n)\}_{n=0}^{\infty}$}
\par
\smallskip
In order to carry out the white noise distribution theory
for the spaces $[\cv]_{\a}$ and $[\cv]_{\a}^{*}$ the
following three conditions have been imposed on
$\{\a(n)\}_{n=0}^{\infty}$ in \cite{kks}:
\par
\smallskip
\begin{enumerate}
\item[(c-1)] There exists a constant $c_{1}$ such that for
any $n\leq m$,
$$\a(n) \leq c_{1}^{m} \a(m).$$
\item[(c-2)] There exists a constant $c_{2}$ such that
for any $n$ and $m$,
$$\a(n+m) \leq c_{2}^{n+m} \a(n) \a(m).$$
\item[(c-3)] There exists a constant $c_{3}$ such that
for any $n$ and $m$,
$$\a(n) \a(m) \leq c_{3}^{n+m} \a(n+m).$$
\end{enumerate}
\par
Note that $c_{i}\geq 1$ for all $i=1, 2, 3$ since
$\a(0)=1$. As shown in Section 3 in \cite{kks}, condition
(c-3) implies condition (c-1). Moreover, it has been proved
in Theorem 4.8 in \cite{kks} that the Bell numbers
$\{b_{k}(n)\}_{n=0}^{\infty}$ satisfy conditions (c-1),
(c-2), and (c-3). Below we give further comments on the
constants $c_{1}, c_{2}$, and $c_{3}$.
\par
Obviously, if a sequence $\{\a(n)\}_{n=0}^{\infty}$ is
non-decreasing, then condition (c-1) is satisfied and
$c_{1}=1$ is the best constant satisfying condition (c-1).
\par
From Eq.~(7.5) in \cite{cks} we have the formula for the
sequence $\{B_{k}(n)\}_{n=0}^{\infty}$ defined in
Eq.~(\ref{eq:n-9}) for $k\geq 2$:
\begin{equation} \label{eq:n-10}
B_{k}(n) = \sum_{j=0}^{\infty} {B_{k-1}(j) \over j!}
\,j^{n},
\end{equation}
where $B_{1}(n)=1$ for all $n$. On the other hand, we can
differentiate both sides of Eq.~(\ref{eq:n-9}) and then
compare the coefficients of $x^{n}$ to get the formula:
\begin{equation} \label{eq:n-11}
B_{k}(n+1) = \sum_{j_{1}+\cdots+j_{k}=n} \> {n! \over
j_{1}! \cdots j_{k}!}\, B_{1}(j_{1}) \cdots B_{k}(j_{k}).
\end{equation}
We see from either Eq.~(\ref{eq:n-10}) or (\ref{eq:n-11})
that the sequence $\{B_{k}(n)\}_{n=0}^{\infty}$ is
increasing. But $b_{k}(n)=B_{k}(n)/\exp_{k}(0)$ and so the
sequence $\{b_{k}(n)\}_{n=0}^{\infty}$ is also increasing.
Hence the Bell numbers satisfy condition (c-1) and the best
constant for $c_{1}$ is $c_{1}=1$.
\par
As mentioned at the end of Section \ref{sec:1}, the Bell
numbers of any order $k\geq 2$ satisfy the inequalities:
$$b_{k}(n) b_{k}(m) \leq b_{k}(n+m) \leq
  2^{n+m} b_{k}(n) b_{k}(m).$$
Hence the Bell numbers satisfy conditions (c-2) and (c-3)
with $c_{2}=2$ and $c_{3}=1$. Obviously, $c_{3}=1$ is the
best constant for condition (c-3). As for the best constant
for $c_{2}$ we have the following
\par
\smallskip\noindent
{\bf Conjecture.} {\em The best constant $c_{2}$ in the
condition {\rm (c-2)} for the Bell numbers
$\{b_{k}(n)\}_{n=0}^{\infty}$ of any order $k\geq 2$ is
$c_{2}=2$.}
\par
\smallskip
Here we prove that the conjecture is true for $k=2$.
It follows from Theorem \ref{thm:3} that $b_{2}(n+m) \leq
2^{n+m} b_{2}(n) b_{2}(m)$. Hence the best constant
$c_{2}$ must be $c_{2}\leq 2$. On the other hand, by
Theorem 4.3 in \cite{kks},
$$\log b_{2}(n) = n\log n - n\log\log n -n + o(n).$$
From this equality we obtain that
$$\log b_{2}(2n) - 2 \log b_{2}(n) = 2n\log 2 -
2n\big(\log\log (2n) - \log\log n\big) + o(n).$$
Then we get the following limit
\begin{equation} \label{eq:n-12}
\lim_{n\to\infty} {1 \over 2n} \log {b_{2}(2n) \over
b_{2}(n)^{2}} = \log 2.
\end{equation}
Now, put $m=n$ in condition (c-2) to get $b_{2}(2n)\leq
c_{2}^{2n} b_{2}(n)^{2}$. This inequality implies that
for all $n\geq 1$,
\begin{equation} \label{eq:n-13}
{1 \over 2n} \log {b_{2}(2n) \over b_{2}(n)^{2}}
\leq \log c_{2}.
\end{equation}
Obviously, Eqs.~(\ref{eq:n-12}) and (\ref{eq:n-13}) show
that $\log 2 \leq \log c_{2}$. Hence $c_{2}\geq 2$. But
we already noted above that $c_{2}\leq 2$. Therefore,
$c_{2}=2$.

\end{document}